\documentclass[aps,pre,twocolumn,a4paper,showkeys,showpacs]{revtex4}%
\usepackage{color}
\usepackage{amsthm}
\usepackage{amsmath}
\usepackage{graphicx}
\usepackage{amssymb}
\usepackage{esint}
\usepackage[unicode=true, pdfusetitle,
bookmarks=true,bookmarksnumbered=false,bookmarksopen=false,
breaklinks=false,pdfborder={0 0 0}%
,backref=false,colorlinks=true,anchorcolor=black,   citecolor=black,
filecolor=black,
menucolor=black,pagecolor=black,urlcolor=black,linkcolor=black]{hyperref}%
\usepackage{amsfonts}
\providecommand{\U}[1]{\protect\rule{.1in}{.1in}}
\pdfoutput=1
\makeatletter
\theoremstyle{plain}
\newtheorem{thm}{Theorem}
\theoremstyle{plain}
\newtheorem{prop}[thm]{Proposition}

\makeatother
\begin{document}
\title{On the Stability of a Chain of Phase Oscillators}
\date{\today}
\author{Jan Sieber}
\affiliation{Department of Mathematics, University of Portsmouth, U.K.}
\author{Tam{\'a}s \surname{Kalm{\'a}r-Nagy}}
\affiliation{Department of Aerospace Engineering, Texas A\&M University}

\begin{abstract}
We study a chain of $N+1$ phase oscillators with asymmetric but  uniform
coupling. This type of chain possesses $2^{N}$ ways to  synchronize in
so-called travelling wave states, i.e. states where  the phases of the single
oscillators are in relative equilibrium. We  show that the number of unstable
dimensions of a travelling wave  equals the number of oscillators with
relative phase close to  $\pi$. This implies that only the relative
equilibrium corresponding  to approximate in-phase synchronization is locally
stable. Despite  the presence of a Lyapunov-type functional periodic or
chaotic phase  slipping occurs. For chains of length $3$ and $4$ we locate
the  region in parameter space where rotations (corresponding to phase
slipping) are present.

\end{abstract}

\maketitle

\section{Introduction}

Investigations of synchronized behavior of coupled nonlinear oscillators
permeate physics \cite{strogatz1988csl,kuramoto1985cdo}, chemistry
\cite{bareli1985scc,crowley1989ets}, biology \cite{kawato1980tcn} and
engineering \cite{ijspeert2001ccp}. Dating back to Huyghens, this problem has
been explored both mathematically and experimentally. The groundbreaking work
of Kuramoto \cite{kuramoto2003cow} led to the realization that synchronization
is a ubiquitous behavior in nature and can be explained by simple models for
interaction between components of a system (see for example the excellent
review of Strogatz \cite{strogatz2000kce}).

Cohen \emph{et al.} \cite{cohen1982ncs} studied a chain of Kuramoto
oscillators to explain the so-called \textit{fictive swimming} observed in the
\textit{Central Pattern Generator} (CPG) of the primitive vertebrate lamprey.
Their purpose was to explain the uniform phase lag along the segmental
oscillators. Later, Kopell and Ermentrout \cite{kopell1988coa} proposed a more
realistic model for fictive swimming and Williams \emph{et al.}
\cite{williams1990fcn} contrasted this model with experimental observations.
Studies of the CPG led to bio-inspired applications in robotics, for example,
the autonomous mobile robotic worm of Conradt and Varshavskaya
\cite{conradt2003dcp} and a turtle-like underwater vehicle by Seo \emph{et
al.} \cite{seo2008cbc}.

Cohen \emph{et al.} \cite{cohen1982ncs} observed that in a chain of
oscillators the oscillators at the end points play as special role:
adapting the natural frequency of the end points controls the phase
shift between neighboring oscillators throughout the entire chain in
the phase-locked equilibrium. (All equilibria of the one-dimensional
homogeneous Kuramoto model with periodic boundary conditions were
found in \cite{mehta2011stationary}.)

In this paper we vary this detuning $\delta$ of the end points and the
coupling strength (more precisely the ratio $k$ between coupling strength down
the chain and up the chain) to study transitions between phase locked
solutions and phase slipping. If one introduces the phase differences between
the oscillators as the new dependent variables then the phase-locked solutions
are equilibria, and for a chain of length $N+1$ there exist $2^{N}$ of these
equilibria. A transition from one phase locked solution to another then
corresponds to a motion along a heteroclinic connection between the
corresponding saddles in the phase space. The $2^{N}$ equilibria can be
classified by an index quantity $\nu_{\pi}$ which counts how many of the phase
differences are equal to $\pi-\delta$. This index turns out to be identical to
the number $\nu_{u}$ of unstable dimensions of the saddle (provided $k>-1$)
such that connections between equilibria of decreasing index $\nu_{\pi}$ are
generic. Solutions with continuously slipping phases between oscillators show
up as rotating waves in the phase space. Periodic phase slipping and its
bifurcations can be computed directly if one assumes, for example, pure phase
coupling (coupling of the type $K\sin(\theta_{j}-\theta_{j-1})$ for
neighboring oscillators $\theta_{j}$ and $\theta_{j+1}$). One of the
co-dimension $1$ boundaries of rotating waves are non-generic connections
between equilibria of identical index $\nu_{\pi}$. We present numerical
evidence of this for the cases of $N=2$ and $N=3$ (that is, the case of $3$
and $4$ oscillators, respectively). A noticeable difference between these
cases is that the basin of attraction for rotating waves appears to be a
slightly smaller fraction of the phase space for larger $N$, and the parameter
region in the $(k,\delta)$-plane permitting rotations (that is, continuous
phase slipping) is smaller for the larger $N$. We conjecture that making an
oscillator chain longer does not make its propensity for phase slipping larger
even if the coupling strength between neighboring oscillators is not increased.

\section{Model Description}

\label{sec:mod}
\begin{figure}[th]
\includegraphics[width=0.48\textwidth]{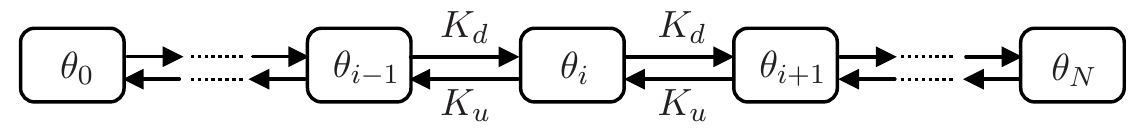}
\caption{Schematic representation of chain coupled oscillators}%
\label{fig:chaintopology}%
\end{figure}

We consider a chain of $N+1$ phase oscillators with nearest
neighbor coupling as shown in Fig.~\ref{fig:chaintopology}.This model, first proposed by Cohen \emph{et al.}
\cite{cohen1982ncs} to explain the 'fictive swimming' observed in Lamprey
spinal cord, is given by
\begin{equation}%
\begin{split}
\frac{d}{dt}\theta_{0}  &  =\omega_{0}+K_{u}\Gamma(\theta_{1}-\theta
_{0})\mbox{,}\\
\frac{d}{dt}\theta_{j}  &  =\omega_{j}+K_{d}\Gamma(\theta_{j-1}-\theta
_{j})+K_{u}\Gamma(\theta_{j+1}-\theta_{j})\\
& \qquad\mbox{\ for \ensuremath{j=1,2,\ldots,N-1},}\\
\frac{d}{dt}\theta_{N}  &  =\omega_{N}+K_{d}\Gamma(\theta_{N-1}-\theta
_{N})\mbox{,}
\end{split}
\label{eq:oscillators}%
\end{equation}
where $\theta_{j}$ and $\omega_{j}\neq0$ are the phase and natural frequency
of the $j$-th oscillator, respectively. The coupling gains are denoted by
$K_{u}$ and $K_{d}$. The prototypical example for the coupling function is
$\Gamma(x)=\sin x$. This choice of the coupling function has the following
essential features:

\begin{enumerate}
\item $\Gamma$ is $2\pi$-periodic, $\left\vert \Gamma(x)\right\vert \leq1$,

\item $\Gamma(x)$ has odd symmetry about $x=0$ (that is, $\Gamma
(x)=-\Gamma(-x)$) and even symmetry about $x=\pi/2$ (that is, $\Gamma
(\pi/2-x)=\Gamma(\pi/2+x)$),

\item $\Gamma^{\prime}(x)>0$ for $x\in\lbrack0,\pi/2]$.
\end{enumerate}



\section{Phase Locking and Traveling Waves}

\label{sec:phaselock}Since the frequencies $\omega_{j}$ are non-zero, system
\eqref{eq:oscillators} does not in general have fixed points. If all
oscillators are identical (that is, $\omega_{0}=\ldots=\omega_{N}$) Equation
\eqref{eq:oscillators} admits a fully synchronised state $\theta_{0}%
=\ldots=\theta_{n}$. Here we focus on the so-called phase locked solutions,
that is, solutions of the form $\theta_{i}-\theta_{j}=\operatorname{const}$.
Of particular importance are phase locked solutions where the phase
differences on the chain are identical, called \textit{uniform traveling
waves.}

To analyze phase locked solutions we introduce the \emph{phase differences}
\begin{equation}
x_{j}=\theta_{j-1}-\theta_{j}%
\mbox{\quad(\ensuremath{j=1\ldots N})},\label{eq:phasediff}%
\end{equation}
as new variables, as well as the rescaled frequency differences $\Omega_{j}$,
the coupling ratio $k$ and the rescaled time $t_{\mathrm{new}}$:
\[
\Omega_{j}=\frac{\omega_{j-1}-\omega_{j}}{K_{u}}%
\mbox{\quad(\ensuremath{j=1\ldots
N}),}k=\frac{K_{d}}{K_{u}}\mbox{,}t_{\mathrm{new}}=K_{u}t\mbox{.}
\]
System \eqref{eq:oscillators} leads to a system of $N$ equations describing
the dynamics of the variables $x_{j}$ wrt the new time $t_{\mathrm{new}}$
\begin{equation}%
\begin{split}
\dot{x}_{1} &  =\Omega_{1}-(1+k)\Gamma(x_{1})+\Gamma(x_{2})\mbox{,}\\
\dot{x}_{i} &  =\Omega_{i}+k\Gamma(x_{i-1})+\Gamma(x_{i+1})-(1+k)\Gamma
(x_{i})\mbox{,}\\
&  \mbox{\qquad\ensuremath{i=2,3,\dots,N-1},}\\
\dot{x}_{N} &  =\Omega_{N}+k\Gamma(x_{N-1})-(1+k)\Gamma(x_{N})\mbox{,}
\end{split}
\label{eq:phase}%
\end{equation}
where $x_{i}(t)$ is on the unit circle $\mathbb{S}^{1}$. The vector form
reads
\begin{equation}
\dot{\mathbf{x}}=\mathbf{\Omega}+C\mathbf{\boldsymbol{\Gamma}}(\mathbf{x}%
)\mbox{,}\label{eq:phasevec}%
\end{equation}
where $\mathbf{x}=[x_{1},x_{2},\cdots,x_{N}]^{T}$, $\mathbf{\Omega}%
=[\Omega_{1},\Omega_{2},\cdots,\Omega_{N}]^{T}$, $\mathbf{\boldsymbol{\Gamma}%
}(\mathbf{x})=[\Gamma(x_{1}),\ldots,\Gamma(x_{N})]$, and $C$ is an $N\times N$
matrix of the form
\begin{equation}
C=%
\begin{bmatrix}
-(1+k) & 1 & \cdots & 0\\
k & -(1+k) & \ddots & 0\\
\vdots & \ddots & -(1+k) & 1\\
0 & 0 & k & -(1+k)
\end{bmatrix}
.\label{eq:-1}%
\end{equation}
The fixed points of Equation (\ref{eq:phasevec}) are given by $\mathbf{\Gamma
}(\mathbf{x})=-C^{-1}\mathbf{\mathbf{\Omega}}$ ($C$ is invertible for
$k\neq-1$, see Proposition~\ref{thm:hyp}) and correspond to the phase locked
solutions satisfying $\theta_{i}-\theta_{j}=\operatorname{const}$ of the
original system (\ref{eq:oscillators}).

If all components of $-C^{-1}\mathbf{\mathbf{\Omega}}$ are less than $1$ in
absolute value, system~\eqref{eq:phasevec} has $2^{N}$ equilibria, because of
the even symmetry about $\pi/2$ ($\Gamma(\pi/2-x)=\Gamma(\pi/2+x)$, see
Section 2) of the coupling function. A uniform traveling wave solution (phase
locked solution with identical phase differences, i.e. $\mathbf{x}%
=[\delta,\cdots,\delta]^{T}$)) exists only if the following conditions
on the rescaled frequency differences $\Omega$ are satisfied:
\begin{equation}
\begin{aligned}\Omega_{1\,\phantom{-1}} & =k\Gamma(\delta)\mbox{,}\\ \Omega_{i\ \phantom{-1}} & =0\mbox{,} & i & =2,\dots,N-1\mbox{,}\\ \Omega_N & =\Gamma(\delta)\mbox{.}\end{aligned}\label{eq:Omcond}%
\end{equation}
This means that the frequencies $\omega_{j}$ of the original system
\eqref{eq:oscillators} must be of the form
\begin{equation}
\begin{aligned}\omega_{0} & =\omega+K_{d}\Gamma(\delta)\mbox{,}\\ \omega_{j} & =\omega\mbox{,} & j & =1,\dots,N-1\mbox{,}\\ \omega_{N} & =\omega-K_{u}\Gamma(\delta)\mbox{,}\end{aligned}\label{eq:omegacond}%
\end{equation}
where $\omega\in\mathbb{R}$. In other words, all oscillators must have
identical natural frequencies except for the two oscillators at the boundary.
The {}``detunings'' (difference from the uniform frequency $\omega$) of the
first and last oscillators are related to one another via the coupling
strengths $K_{u}$ and $K_{d}$.

The two primary parameters affecting the dynamics are the coupling strength
ratio $k$ and $\delta$. Without loss of generality we can restrict our
considerations to the parameter set

\begin{itemize}
\item $\delta\in[0,\pi/2]$: only $\Gamma(\delta)$ enters the equation  and for
negative $\delta$ we can apply the transformation  $x\mapsto-x$ since $\Gamma$
is odd.

\item $k\geq-1$: for $k<-1$ we can apply the transformation  $[x_{1}%
,\ldots,x_{N}]_{\mathrm{new}}=[x_{N},\ldots,x_{1}]_{\mathrm{old}}$,
$t_{\mathrm{new}}=k_{\mathrm{old}}t$, and  $k_{\mathrm{new}}%
=-1/k_{\mathrm{old}}$ such that  $k_{\mathrm{new}}>0$. Note that this
transformation involves  reversal of time direction (since $k_{\mathrm{old}}$
is negative) so  statements on stability will have to be replaced by the
corresponding statements on instability and vica versa.
\end{itemize}

We observe that the coupling \emph{strength} does not enter
Equation~\eqref{eq:phase}, which determines the dynamics, at all. Only the
\emph{ratio} between down-chain and up-chain coupling strength matters. The
absolute value of the coupling strength then determines the time scale with
respect to the the original time (Equation~\eqref{eq:phase} is with respect to
a rescaled time $K_{u}t$).

\section{Stability of Traveling Waves}

\label{sec:mi}Naturally, we are interested in characterizing the local
stability of traveling waves, that is, of equilibria of \eqref{eq:phasevec}.
We note that even though the uniform traveling waves can occur only in a
slightly degenerate parameter setting, these waves are robust in the sense
that slightly non-uniform traveling waves will exist for slight perturbations
of these parameters values. For example, non-uniform traveling waves exist for
a linear gradient frequency distribution ($\Omega=\operatorname{const}.$).
Cohen \emph{et al.} \cite{cohen1982ncs} derived a necessary condition for the
existence of the phase-locked solutions. Ermentrout and Kopell
\cite{ermentrout1984fpc} showed the existence of \emph{frequency plateaus}
when this necessary condition is violated.

The linearization of \eqref{eq:phasevec} about an equilibrium $x_{*}$ is
$J=C\mathbf{\boldsymbol{\Gamma\mathrm{^{\prime}}}}(x_{*})$ where $C$ is the
coupling matrix given in \eqref{eq:-1}, and%
\begin{equation}
\boldsymbol{\Gamma}^{\prime}(x_{*})=\rho\cdot\operatorname{diag}(\sigma
_{1},\dotsc,\sigma_{n}),
\end{equation}
is a diagonal matrix with $\rho=\Gamma^{\prime}(\delta)>0$ (since $\delta
\in[0,\pi/2]$) and $\sigma_{i}=\pm1$, depending on whether the $i$th component
of $x_{*}$ is equal to $\delta$ or $\pi-\delta$. Even though the eigenvalues
of the matrix $C$ are known, i.e.%

\begin{equation}
\lambda_{j}=-(1+k)+2\sqrt{k}\cos\frac{j\pi}{N+1},\mbox{}j=1,\dots,N,
\end{equation}
the eigenvalues of the product $C\mathbf{\boldsymbol{\Gamma\mathrm{^{\prime}}%
}}(x_{\ast})$ are not known analytically. Therefore, the goal of this Section
is to establish a simple connection between the number of stable and unstable
directions of an equilibrium $x_{\ast}$ (and thus its stability) of
\eqref{eq:phasevec} and the number of components of $x_{\ast}$ equal to
$\delta$ and $\pi-\delta$. To establish the result, we associate the four
indices with each equilibrium $x_{\ast}$
\begin{equation}%
\begin{split}
\nu_{0}(x_{\ast}) &
=\ \mbox{number of components of \ensuremath{x_{*}} equal to \ensuremath{\delta},}\\
\nu_{\pi}(x_{\ast}) &
=\ \mbox{number of components of \ensuremath{x_{*}} equal to \ensuremath{\pi-\delta},}\\
\nu_{u}(x_{\ast}) &
=\ \mbox{dimension of unstable manifold of \ensuremath{x_{*}},}\\
\nu_{s}(x_{\ast}) &
=\ \mbox{dimension of stable manifold of \ensuremath{x_{*}}.}
\end{split}
\label{eq:class}%
\end{equation}
The relation between these indices is given by the following

\begin{thm}
[Dimension of invariant subspaces]\label{thm:invdim} Let $x_{*}$ be an
equilibrium of \eqref{eq:phasevec}. If $k>-1$ then
\begin{equation}
\nu_{0}(x_{*})=\nu_{s}(x_{*})\mbox{\ and\ }\nu_{\pi}(x_{*})=\nu_{u}%
(x_{*})\mbox{.}\label{eq:kgtm}%
\end{equation}

\end{thm}

\paragraph*{Proof}

Relation \eqref{eq:kgtm} is proven indirectly via the Lyapunov-type
functional
\begin{equation}
E(x)=\sum_{i=1}^{N}\left[  \int_{0}^{x_{i}}\Gamma(y)\,y-\Gamma(\delta
)\,x_{i}\right]  \mbox{.}\label{eq:efunc}%
\end{equation}
We observe that the time derivative of $E$ along trajectories of
\eqref{eq:phasevec} is
\begin{equation}%
\begin{split}
\dot{E}= &  \sum_{i=1}^{N}\left[  \Gamma(x_{i})-\Gamma(\delta)\right]  \dot
{x}_{i}\\
= &  -\frac{k+1}{2}\left(  \left(  \Gamma(x_{1})-\Gamma(\delta)\right)
^{2}+\left(  \Gamma(x_{N})-\Gamma(\delta)\right)  ^{2}%
+\phantom{\sum_i^N}\right.  \\
&  \left.  +\sum_{i=1}^{N-1}\left(  \Gamma(x_{i})-\Gamma(x_{i+1})\right)
^{2}\right)  .
\end{split}
\label{eq:dedt}%
\end{equation}
If $k=-1$ then $E(x(t))$ is constant along trajectories, however it is
strictly decreasing for $k>-1$ provided $x(t)$ is not an equilibrium of
\eqref{eq:phasevec} (i.e. whenever at least one of the phase differences
$x_{i}$ does not satisfy $\Gamma(x_{i})=\Gamma(\delta)$). With this functional
$E$ and a bound on its rate of decrease \eqref{eq:dedt} we establish in
Appendix \ref{sec:smallproofs} that the linearization
$J=C\mathbf{\boldsymbol{\Gamma\mathrm{^{\prime}}}}(x_{\ast})$ is hyperbolic
for $k>-1$. This implies that eigenvalues of $J$ cannot \emph{cross }the
imaginary axis when $k$ is varied. Thus, for all $k>-1$ the indices $\nu
_{u}(x_{\ast})$ and $\nu_{s}(x_{\ast})$ are the same as for $k=0$. Since for
$k=0$ the Jacobian $J$ is an upper diagonal matrix with its eigenvalues on its
diagonal, the number of negative diagonal entries of $J$ is $\nu_{0}(x_{\ast
})$ (the diagonal entry of $C$ is $-1$), and the number of positive diagonal
entries of $J$ is $\nu_{\pi}(x_{\ast})$. \hfill$\square$

We observe that for $\delta=0$ the functional $E(x)$ is bounded, leading for
$k\neq1$ to the result that full synchronization (the equilibrium
$(0,\ldots,0)$ for $\delta=0$) is globally stable in the same way as the
classical results \cite{hoppensteadt1997weakly}.

\section{Rotating waves and saddle connections}

\label{sec:rot}Despite the existence of a Lyapunov functional the dynamics of
\eqref{eq:phasevec} is not necessarily trivial because the phase space is an
$N$-dimensional torus. Generically we can expect that the unstable manifold of
an equilibrium $x$ and the stable manifold of an equilibrium $y$ intersect,
giving rise to heteroclinic saddle connections, if $\nu_{\pi}(x)+\nu_{0}%
(y)>N$. This implies that for any equilibrium $x$ heteroclinic connections to
all equilibria $y$ satisfying $\nu_{\pi}(y)<\nu_{\pi}(x)$ are generic.

Saddle connections between equilibria of the same type, that is $\nu_{\pi
}(x)=\nu_{\pi}(y)$, are of co-dimension $1$ and can be achieved by tuning the
parameters $\delta$ and $k$. These non-generic saddle connections also form
co-dimension $1$ boundaries of periodic orbits of rotating wave type.
\begin{figure}[ptb]
\includegraphics[width=0.48\textwidth]{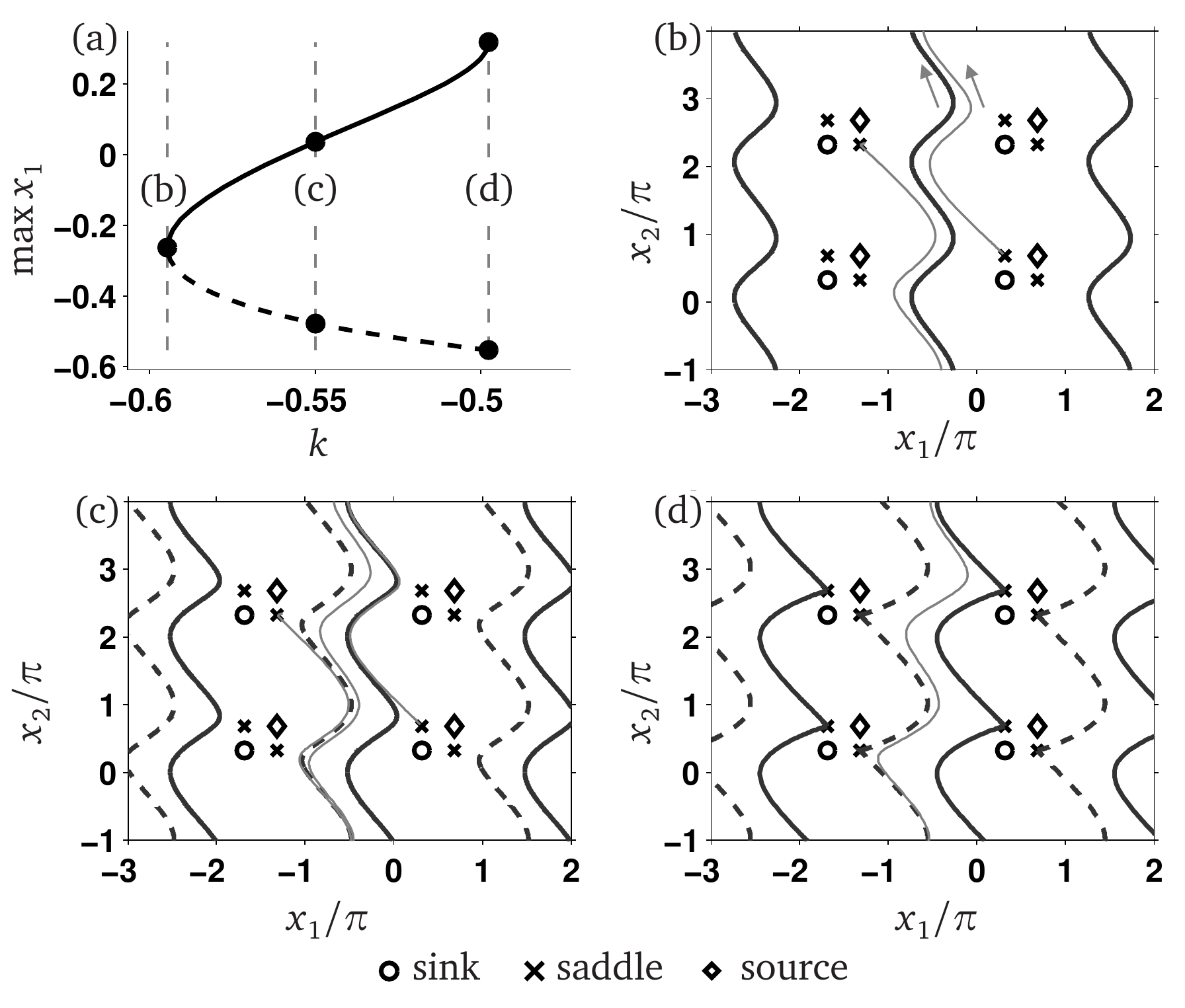}\caption{Co-dimension $1$
homoclinic connections and rotating waves for $N=2$, $\Gamma(x)=\sin x$ and
$\delta=1$. (a):~bifurcation diagram for rotating wave solutions (computed
with (r)AUTO \cite{DCFKSW98,S07}), (b)--(d):~phase portraits showing
equilibria homoclinic connections and rotating waves for parameter values of
$k$ indicated in panel (a). The direction of the flow is always upward.}%
\label{fig:bifpp2}%
\end{figure}

Figure~\ref{fig:bifpp2}(a) shows for $N=2$ that periodic orbits (of rotating
wave type) are possible for a range of parameters $k$ and $\delta=1$ despite
the existence of a Lyapunov functional $E$. In the parameter space the
existence of periodic orbits is bounded by two co-dimension $1$ bifurcations:
a saddle-node of periodic orbits corresponding to the phase portrait in
Figure~\ref{fig:bifpp2}(b) at the minimal $k$ supporting periodic orbits and a
heteroclinic connection between neighboring saddles (the $x_{2}$-component has
increased by $2\pi$ along the heteroclinic connection). \begin{figure}[ptb]
\includegraphics[width=0.4\textwidth]{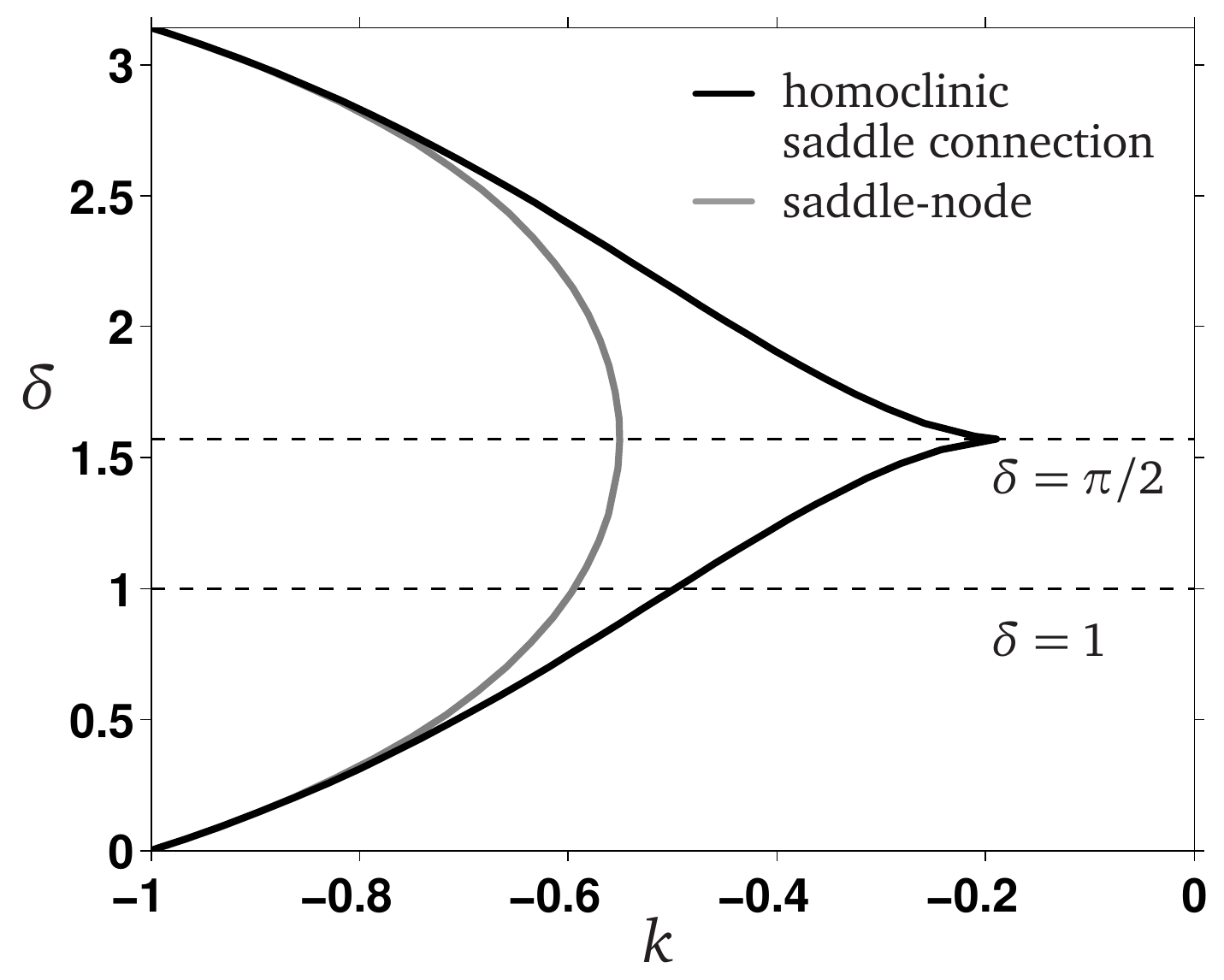}\caption{Bifurcation diagram
in the $(k,\delta)$-plane for $N=2$ and $\Gamma(x)=\sin x$ (computed with
(r)AUTO \cite{DCFKSW98,S07}).}%
\label{fig:bif2}%
\end{figure}

If the phase shift by $2\pi$ is ignored then the saddle connection is
homoclinic. Figure~\ref{fig:bif2} shows the parameter region in the
$(k,\delta)$-plane where periodic orbits of rotating wave type exist. The
special points in the plane are at $\delta=\pi/2$ when all four equilibria
collapse into a single degenerate equilibrium and at $k=-1$ when the system is
a single-degree-of-freedom conservative oscillator for $N=2$. We also note
that the diagram is reflection symmetric about the line $\delta=\pi/2$ as expected.

\begin{figure}[ptb]
\includegraphics[width=0.4\textwidth]{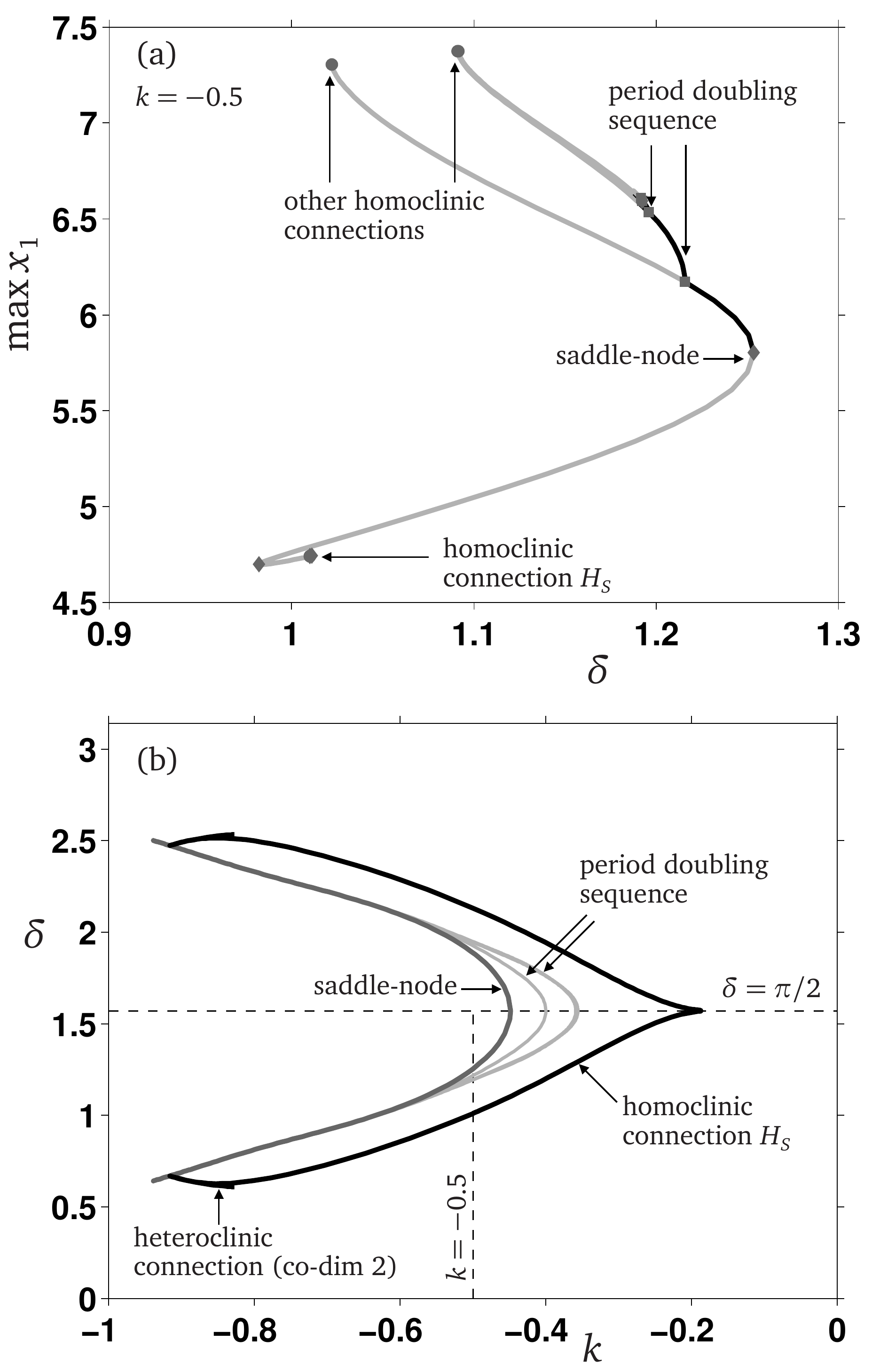}\caption{Bifurcation diagrams
for $N=3$ and $\Gamma(x)=\sin x$ (computed with (r)AUTO \cite{DCFKSW98,S07}).
Panel~(a): family of periodic orbits for fixed $k=-0.5$ and varying $\delta$
and its codimension-$1$ bifurcations. Panel~(b): region with rotating waves in
the $(k,\delta)$-plane with its codimension-$1$ bifurcations. }%
\label{fig:bif3}%
\end{figure}

In dimensions $N>2$ the dynamics can become chaotic. Numerical evidence for
this is shown in Figure~\ref{fig:bif3} for $N=3$. Panel (a) shows a family of
periodic orbits of rotating wave type $(x_{1}(T),x_{2}(T),x_{3}(T))=(x_{1}%
(0),x_{2}(0)+2\pi,x_{3}(0)+2\pi)$ where $T$ is the period, $\delta$ is varied,
and $k=-0.5$ is kept fixed. We observe that this family undergoes a period
doubling sequence. Moreover, the lower (predominantly unstable) part of the
branch reaches a homoclinic connection to the saddle $x_{s}=(\pi-\delta
,\delta,\delta)$ at $H_{s}$. The eigenvalues at this saddle have the form
$(\mu_{+},-\mu_{-}\pm i\omega_{-})$ where $\mu_{\pm}>0$ and $\omega_{-}>0$ are
real numbers and $\mu_{-}<\mu_{+}<2\mu_{-}$. Shil'nikov's results imply that
there is an infinite number of period doubling cascades of stable periodic
orbits close to the homoclinic connection $H_{s}$ (under certain
non-degeneracy conditions, see \cite{SSTC01,K04}). The precise sequence of
period-$n$ branches for a homoclinic to a saddle of this type, and how to
calculate them numerically can be found in \cite{OKC00}. The other end of the
family of periodic orbits is also a homoclinic connection, to the saddle
$(\delta,\pi-\delta,\delta)$, which has three real eigenvalues $\mu_{+}%
>0>-\mu_{-,1}>-\mu_{-,2}$ where $\mu_{+}>\mu_{-,1}$. This implies that there
is only a single periodic orbit close to this homoclinic connection (no
\emph{snaking}), and that this periodic orbit is unstable \cite{K04}.

\begin{figure}[ptb]
\includegraphics[width=0.48\textwidth]{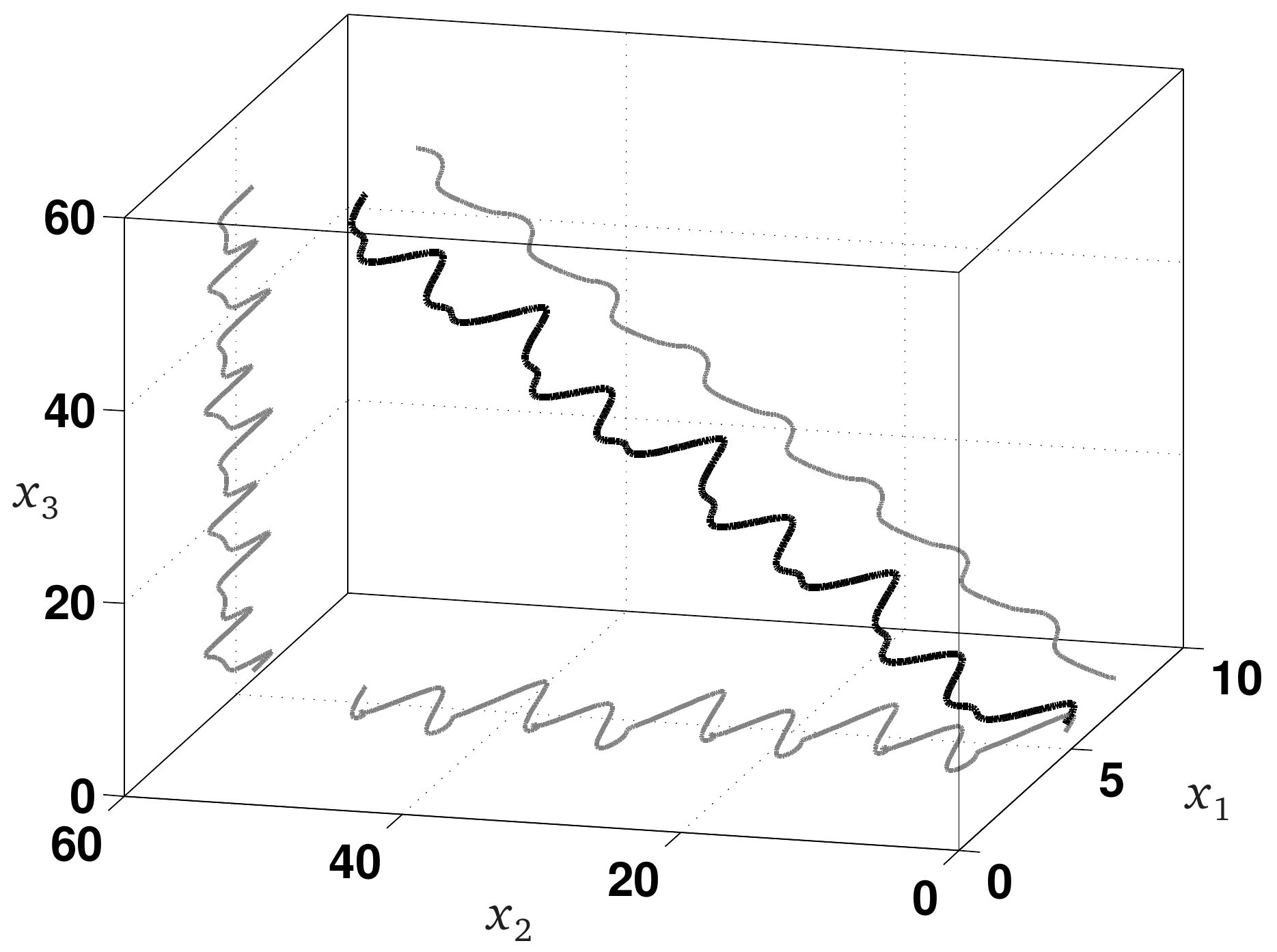}\caption{Rotating periodic
orbit of period $8$ for $N=3$, $\Gamma(x)=\sin x$, $k=-0.5$, $\delta
\approx1.19$, near the end of the period doubling sequence in parameter space.
}%
\label{fig:phas3}%
\end{figure}

In the two-parameter plane (shown in Figure~\ref{fig:bif3}(b)) we observe a
shape that looks superficially similar to the case $N=2$ except that the range
of $\delta$ over which the region of rotating waves extend is smaller.
However, as one can see from Figure~\ref{fig:bif3}(a) the rotating waves are
not stable inside the region bounded by the homoclinic connection $H_{s}$ and
the saddle-node. The period doubling sequence (also shown in
Figure~\ref{fig:bif3}(a)) gives a better estimate for the region of stable
periodic rotating waves. The remainder of the region is not filled with
chaotic rotations because the chaotic attractor at the end of the primary
period doubling cascade collides with heteroclinic saddle connections.
Figure~\ref{fig:phas3} shows a periodic motion of period $8$, which is close
to the end of the period doubling sequence in parameter space.

Special points in the bifurcation diagram~\ref{fig:bif3}(b) are on the
symmetry line $\delta=\pi/2$ (where all equilibria collapse to a single
degenerate equilibrium) and various degeneracies of the homoclinic connection,
for example, the heteroclinic connection $(\pi-\delta,\delta,\delta)\to
(2\pi+\delta,\pi-\delta,\delta)\to(\pi-\delta,2\pi+\delta,2\pi+\delta)$
indicated in Figure~\ref{fig:bif3}(b) (this degeneracy ends the numerical
continuation of the homoclinic).

In the original variables $\theta_{j}$ the rotating waves shown in this
Section correspond to a continually drifting phase difference between the
oscillators $2$ and $3$ (for $N=2$) while the phase difference between
oscillators $1$ and $2$ remains bounded. The homoclinic boundary corresponds
to regimes where two neighboring oscillators hover in anti-phase for a long
time before a phase slip occurs. For $N=3$ the chaotic regimes near the
Shil'nikov saddle correspond to the regime where oscillators $1$ and $2$ are
nearly in anti-phase and $2$, $3$ and $4$ are nearly in-phase for a long time
before the phase differences between $2$ and $3$, and $3$ and $4$ slip
(nearly) simultaneously.

\begin{figure}[th]
\centering
\includegraphics[width=0.7\columnwidth]{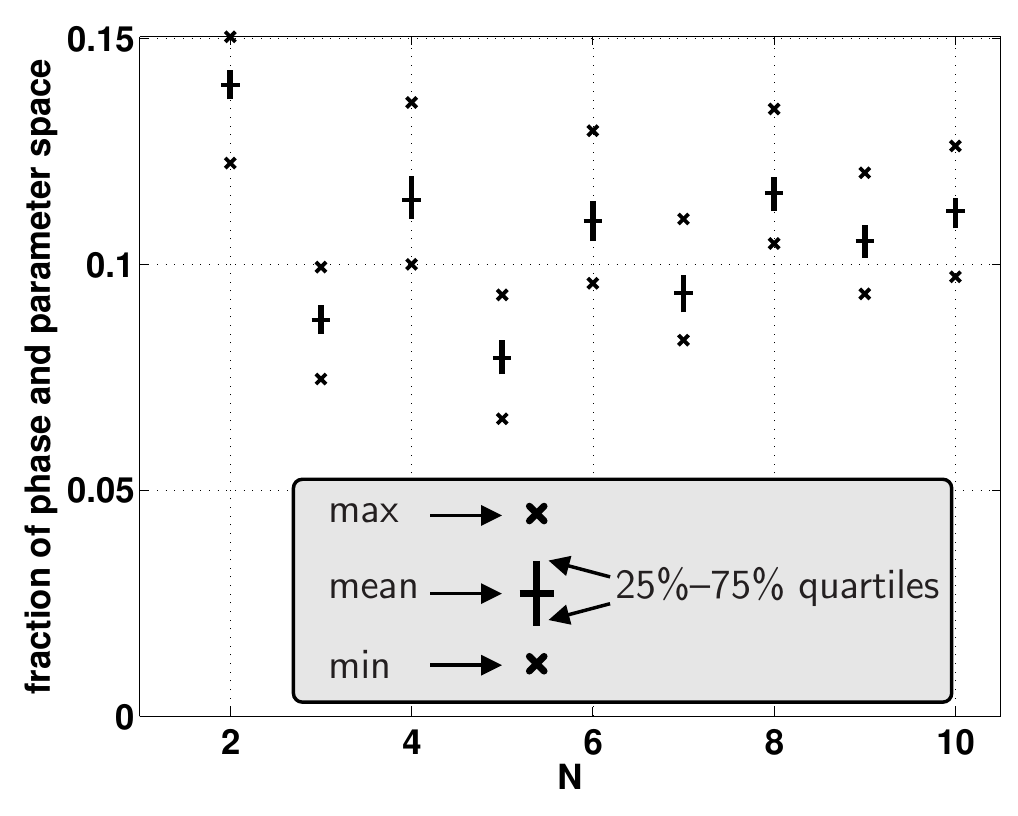}  \caption{Fraction of
phase and parameter space not attracted to synchronization. These are
simulation results obtained during $100$ trials where in each trial we chose
$5,000$ random initial conditions $x(0)\in\mathbb{S}^{N}$ and random
parameters $\delta\in\lbrack0,\pi/2]$, $k\in\lbrack-1,0]$ for small $N\geq2$.
During each trial simulation we record when a trajectory enters the
neighborhood of the stable equilibrium $(\delta,\ldots,\delta)$. We stop the
simulation if no new trajectory enters this neighborhood during a time period
of $200$. The plot (horizontal line) shows which fraction of the $5,000$
realizations was not yet trapped on average when we stopped the simulation.
The vertical lines give the variability between trials ($25\%$ and $75\%$
quartile) and the crosses are the extreme (maximal and minimal) result of the
$100$ trials.}%
\label{fig:slipstats}%
\end{figure}As the number of oscillators is increased by $1$, we observe that
the parameter region which supports stable periodic rotations shrinks. Put
another way, full synchronization becomes more prevalent for $N=3$. The
simulation results shown in Figure~\ref{fig:slipstats} support this
observation. This is in contrast to the common feature that longer chains of
oscillators require stronger coupling for synchronization
\cite{hoppensteadt1997weakly,Hong2005} in symmetrically coupled chains with
random frequencies. The reason behind this apparent difference is that in our
setup the only detuned oscillators are at the boundary (indices $0$ and $N$),
and that the effect of the boundary diminishes for increasing $N$. We did not
observe, however, that the fraction of slipping realizations approaches zero
for large $N$.

\section{Conclusion}

\label{sec:conc}In this paper we studied the dynamical system describing the
phase differences of a uniform chain of oscillators. Two important parameters
(other than the number of elements in the chain) are the detuning of the end
points, $\delta$ and the ratio $k$ of the coupling strengths upward and
downward the chain. We have classified the phase locked states by relating the
number of phase differences that equal $\pi-\delta$ to the number of unstable
directions that this phase locked state has as a saddle in the phase space. We
also studied periodic phase slipping (rotating waves) systematically for $N=2$
($3$ oscillators) and $N=3$ ($4$ oscillators). The most curious difference
between the case of $N=2$ and $N=3$ is the shrinking of the parameter region
where periodic rotations occur. Simulations suggest that their basin of
attraction also shrinks as a fraction of the whole phase space volume. This
would suggest that longer chains have more robust phase locking. This effect
appears to hold with respect to both aspects of robustness combined: the
fraction of initial values (treating parameters also as dynamic variables with
trivial dynamics, $\dot{\delta}=0$, $\dot{k}=0$) in the phase space leading to
locking is smallest for $N=2$.

\begin{acknowledgments}
The authors would like to thank Siming Zhao for his help in this  work. This
material is based upon work supported by the National  Science Foundation
under Grant No. 0846783 and by the US Air Force  Office of Scientific Research
under Grant No. FA9550-08-1-0333.
\end{acknowledgments}

\bibliographystyle{unsrt}
\bibliography{CPG}

\appendix

\section{Proof of hyperbolicity of equilibria}

\label{sec:smallproofs}

\begin{prop}
[Absence of eigenvalues with zero real part]\label{thm:hyp} Let $k>-1$, and
let $x_{*}$ be an equilibrium of \eqref{eq:phasevec}. Then the Jacobian
$J=C\boldsymbol{\Gamma}^{\prime}(x_{*})$ of \eqref{eq:phasevec} in $x_{*}$
cannot have eigenvalues on the imaginary axis.
\end{prop}

\paragraph*{Proof}

First we note that the matrix $J$ is regular for $k>-1$ (and, thus, cannot
have an eigenvalue $0$):
\begin{equation}
\begin{aligned} \det\,&C&=\ &\frac{\left(-1\right)^{N}\left(k^{N+1}-1\right)}{k-1} &\neq&0\mbox{,}\\ \det\,& \boldsymbol{\Gamma}'(x_{*})&=\ &\Gamma'(\delta)^{N}(-1)^{\nu_{\pi}} &\neq&0\mbox{.} \end{aligned}
\end{equation}
We show the absence of purely imaginary eigenvalues indirectly. In preparation
for this part of the proof we establish a bound on the gradient of Lyapunov
functional $E$, given in \eqref{eq:efunc}. In an equilibrium $x_{*}$ the
gradient of $E$ vanishes and the Hessian $H_{*}$ of $E$ in $x_{*}$ is equal to
$\boldsymbol{\Gamma}^{\prime}(x_{*})$:
\begin{equation}%
\begin{split}
H_{*}=H(x_{*}) & =\frac{\partial^{2}}{(\partial x)^{2}}E(x)\vert_{x=x_{*}}\\
& = \rho\cdot\operatorname{diag}(\sigma_{1},\dotsc,\sigma_{n}) =
\boldsymbol{\Gamma}^{\prime}(x_{*}).\label{eq:d2e}%
\end{split}
\end{equation}
Consequently, the equilibrium $(\delta,\dotsc,\delta)$ is a local minimum of
$E$, and $(\pi-\delta,\dotsc,\pi-\delta)$ is a local maximum of $E$.
Equilibria which have some components equal to $\delta$ and others equal to
$\pi-\delta$ are saddle points of the graph of $E$. As we want to study the
local stability of an equilibrium $x_{*}$ we introduce quantities measuring
the deviation from $x_{*}$:
\begin{equation}
y=x-x_{*}\mbox{\ and\
}D(y)=E(y+x_{*})-E(x_{*}),
\end{equation}
such that $D(0)=0$ and
\begin{equation}
D(y)=y^{T}H_{*}y+O(\|y\|^{3}),\label{eq:ebound}%
\end{equation}
for all $y$ in a neighborhood of $0$. Furthermore, equation \eqref{eq:dedt}
estimates $\dot{D}$ from above for $k>-1$:
\begin{equation}
\dot{D}(y)\leq-c_{0}\|y\|^{2}\mbox{,}\label{eq:dedtbound}%
\end{equation}
where $c_{0}>0$ is a constant independent of $y$.

Assume that $J$ has a purely imaginary eigenvalue, say, $\lambda=i\mu$ where
$\mu>0$, and let $u$ be the corresponding eigenvector. We choose the time
$T=2\pi/\mu$. Let $\epsilon>0$ be sufficiently small ($\epsilon$ will depend
on $T$). After time $T$ the solution of \eqref{eq:phasevec} starting from
$x(0)=x_{*}+\epsilon u$ (let us call the solution $x(\cdot)$) satisfies
$x(T)=x_{*}+\epsilon u+O(\epsilon^{2})$. Consequently,
\begin{align}
D(x(T)) & -D(x(0))=\nonumber\\
=\  &  [x(T)-x_{*}]^{T}H_{*}[x(T)-x_{*}]+O(\|x(T)-x_{*}\|^{3})\nonumber\\
&  -\epsilon^{2}u^{T}H_{*}u-O(\epsilon^{3})\nonumber\\
=\  &  O(\epsilon^{3})\mbox{.}\label{eq:dt0}%
\end{align}
On the other hand, the trajectory $x(t)$ for $t\in[0,T]$ is a perturbation of
order $\epsilon^{2}$ of an ellipse with a minimal radius of order $\epsilon$
around $x_{*}$. This means that we can choose a uniform constant $c_{1}$ of
order $1$ such that $c_{1}\epsilon$ is smaller than this minimal radius
($c_{1}$ is uniform in $\epsilon$) and, hence,
\begin{equation}
\|x(t)-x_{*}\|\geq c_{1}\epsilon\label{eq:xta}%
\end{equation}
for all $t\in[0,T]$. Inequality \eqref{eq:dedtbound} implies that
\begin{align*}
D(x(T)) & -D(x(0)) =\int_{0}^{T}\dot{D}(x(t))dt & \\
&  \leq-c_{0}\int_{0}^{T}\|x(t)-x_{*}\|^{2}\, dt &  &
\mbox{(due to \eqref{eq:dedtbound}),}\\
&  \leq-c_{0}Tc_{1}^{2}\epsilon^{2} &  &  \mbox{(due to \eqref{eq:xta}),}
\end{align*}
which contradicts \eqref{eq:dt0} because the constants $c_{0}$ and $c_{1}$ are
independent of $\epsilon$. \hfill$\square$

\end{document}